\documentclass[11pt,reqno,letterpaper]{amsart}

\usepackage{amssymb}
\usepackage{amsmath}
\usepackage{amsthm}
\usepackage{bbm}
\usepackage{doi}
\usepackage{graphicx}

\usepackage{bookmark}
\usepackage{hyperref}
\hypersetup{pdfstartview={FitH}}

\numberwithin{equation}{section}

\newtheorem{theorem}{Theorem}
\newtheorem{lemma}[theorem]{Lemma}

\theoremstyle{definition}
\newtheorem{remark}[theorem]{Remark}

\renewcommand{\leq}{\leqslant}
\renewcommand{\geq}{\geqslant}

\newcommand{\N}{\mathbb{N}}

\begin{document}

\title{On the irrationality of certain super-polynomially decaying series}

\author[T. Crmari\'{c}]{Ton\'{c}i Crmari\'{c}}
\address{T. C., Department of Mathematics, Faculty of Science, University of Split, Ru\dj{}era Bo\v{s}kovi\'{c}a 33, 21000 Split, Croatia}
\email{tcrmaric@pmfst.hr}

\author[V. Kova\v{c}]{Vjekoslav Kova\v{c}}
\address{V. K., Department of Mathematics, Faculty of Science, University of Zagreb, Bijeni\v{c}ka cesta 30, 10000 Zagreb, Croatia}
\email{vjekovac@math.hr}



\begin{abstract}
We give a negative answer to a question by Paul Erd\H{o}s and Ronald Graham on whether the series
\[ \sum_{n=1}^{\infty} \frac{1}{(n+1)(n+2)\cdots(n+f(n))} \]
has an irrational sum whenever $(f(n))_{n=1}^{\infty}$ is a sequence of positive integers converging to infinity. To achieve this, we generalize a classical observation of S\={o}ichi Kakeya on the set of all subsums of a convergent positive series. We also discuss why the same problem is likely difficult when $(f(n))_{n=1}^{\infty}$ is additionally assumed to be increasing.
\end{abstract}

\maketitle


\section{Introduction}

Paul Erd\H{o}s and Ronald Graham commented in their well-known problem book on combinatorial number theory \cite{EG80} at the very end of chapter \emph{Irrationality and Transcendence}:
\begin{quote}
\emph{It seems that series like
\begin{equation}\label{eq:specialseries}
\sum_{n=1}^{\infty} \Big( \prod_{i=1}^{n} (n+i) \Big)^{-1}
\end{equation}
are very hard to treat, though they surely are irrational.} \cite[p.~65]{EG80}
\end{quote}
The series \eqref{eq:specialseries} is quite special as Mathematica \cite{Mathematica} computes its sum 
\[ 0.592296536469326575660415\ldots \]
in terms of the integral resembling the Gauss error function as\footnote{This is also easy to prove directly, by expanding the exponential function in $e^{1/4-t^2}$ into the power series and then observing that the integral $\int_{0}^{1/2} (1/4-t^2)^{n-1}/(n-1)! \,\textup{d}t$ evaluates as the $n$-th term of the series \eqref{eq:specialseries}.}
\[ e^{1/4} \int_0^{1/2} e^{-t^2} \,\textup{d}t, \]
but the authors still do not know of an argument showing that this is, in fact, an irrational number.

Erd\H{o}s and Graham then posed the following general question:
\begin{quote}
\emph{Let $f(n)\to\infty$ as $n\to\infty$. Is it true that
\begin{equation}\label{eq:generalseries}
\sum_{n=1}^{\infty} \Big( \prod_{i=1}^{f(n)} (n+i) \Big)^{-1}
\end{equation}
is irrational?} \cite[p.~66]{EG80}
\end{quote}
Recently this question also appeared as Problem \#270 on Thomas Bloom's website \emph{Erd\H{o}s problems} \cite{EP}.

The reasoning behind the heuristics for the aforementioned question is that the series converges very rapidly starting with the terms for which $f(n)$ becomes very large. Since the series with terms that decay extremely rapidly (e.g., faster than double exponentially) are known to have irrational sums \cite{ES64,E75}, it is also reasonable to suspect that \eqref{eq:generalseries} is never a rational number. However, our main result shows that precisely the opposite is true: the series can sum up to any positive number, rational or irrational.

Let $\N$ denote the set of positive integers.

\begin{theorem}\label{thm:main}
The set
\begin{equation}\label{eq:thewholeset} 
\bigg\{ \sum_{n=1}^{\infty} \frac{1}{\prod_{i=1}^{f(n)}(n+i)} \,:\, (f(n))_{n=1}^{\infty}\in\N^{\N},\ \lim_{n\to\infty} f(n) = \infty \bigg\} 
\end{equation}
is equal to the whole interval $(0,\infty)$.
\end{theorem}

\begin{figure}
\includegraphics[width=0.6\linewidth]{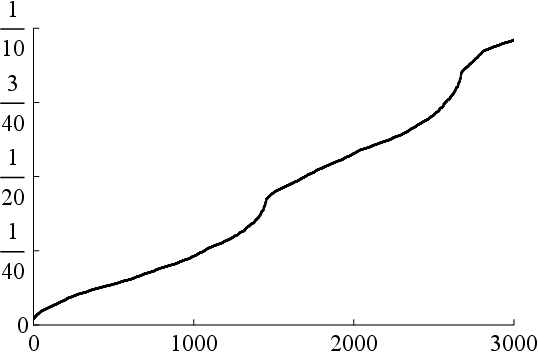}
\caption{Partial series of \eqref{eq:generalseries} with $1\leq f(1)\leq 5$, $2\leq f(2),f(3)\leq 5$, $3\leq f(4),f(5),f(6),f(7)\leq 5$; the first $3000$ sums sorted in the ascending order.}
\label{fig:All}
\end{figure}

In particular, the above general question of Erd\H{o}s and Graham has a negative answer.
It is impossible to depict all sums used in the actual proof, but Figure \ref{fig:All} can already illustrate how the set \eqref{eq:thewholeset} starts filling the whole interval $(0,\infty)$ on the vertical axis.

Interestingly, the proof of Theorem \ref{thm:main} will only require variants of classical ideas on subseries sums, dating back to S\={o}ichi Kakeya \cite{Kak14a,Kak14b} more than a hundred years ago. An extensive reference for similar analytical tools applied to irrationality problems is a recent paper by Terence Tao and one of the present authors \cite{KT25}.

Finally, Erd\H{o}s and Graham also gave a more reserved final comment:
\begin{quote}
\emph{The answer is almost surely in the affirmative if $f(n)$ is assumed to be nondecreasing but at present we lack methods to decide such questions.} \cite[p.~66]{EG80}
\end{quote}
Simple ideas used in the proof of Theorem \ref{thm:main} do not seem to be sufficient to also resolve the latter variant of their question. On the contrary, we prove that the corresponding set of all possible values of \eqref{eq:generalseries} is both measure-theoretically and topologically very different.

\begin{theorem}\label{thm:cantor}
The set
\begin{equation}\label{eq:increasingset}
\bigg\{ \sum_{n=1}^{\infty} \frac{1}{\prod_{i=1}^{f(n)}(n+i)} \,:\, (f(n))_{n=1}^{\infty}\in\N^{\N}\text{ is increasing},\ \lim_{n\to\infty} f(n) = \infty \bigg\} 
\end{equation}
has zero Lebesgue measure and, consequently also, empty interior.
\end{theorem}

\begin{figure}
\includegraphics[width=0.6\linewidth]{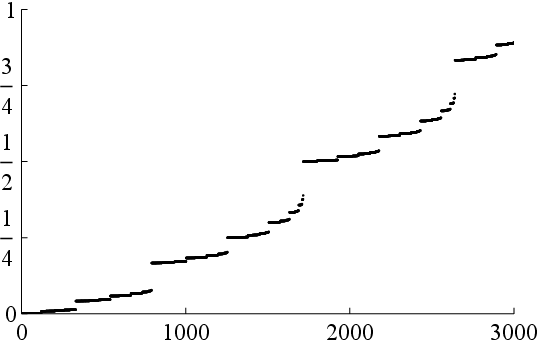}
\caption{Partial series of \eqref{eq:generalseries} with $1\leq f(1)\leq f(2)\leq f(3)\leq f(4)\leq f(5)\leq f(6)\leq f(7)\leq 8$; the first $3000$ sums sorted in the ascending order.}
\label{fig:Inc}
\end{figure}

It is again impossible to depict the whole set \eqref{eq:increasingset}, but Figure \ref{fig:Inc} approximates it by vertical coordinates of $3000$ points and illustrates that it is fractal and ``full of holes.''
This time one should no longer expect an easy negative answer to the irrationality problem, since finding a rational number inside a negligible set is a difficult task.

\section{Preliminary observations}

The following simple result was already known to Kakeya \cite{Kak14a,Kak14b}.

\begin{lemma}[Kakeya]\label{lm:Kakeya}
Let $x=(x_n)_{n=1}^{\infty}$ be a sequence of strictly positive real numbers such that the series $\sum_{n}x_n$ converges.
\begin{itemize}
\item[(a)] If
\begin{equation}\label{eq:Kak1}
\sum_{k=n+1}^{\infty} x_k \geq x_n
\end{equation}
holds for every sufficiently large $n\in\N$, then
\begin{equation}\label{eq:achievement}
\Big\{ \sum_{n=1}^{\infty} \varepsilon_n x_n \,:\, \varepsilon_n\in\{0,1\} \text{ for every } n\in\N \Big\}
\end{equation}
is a finite union of nondegenerate bounded closed intervals. 
Moreover, if condition \eqref{eq:Kak1} holds for every $n\in\N$, then the set \eqref{eq:achievement} is a single interval, namely $[0,\sum_{n=1}^{\infty}x_n]$.
\item[(b)] If
\begin{equation}\label{eq:Kak2}
\sum_{k=n+1}^{\infty} x_k < x_n
\end{equation}
holds for every sufficiently large $n\in\N$, then \eqref{eq:achievement} is a closed set with empty interior. In fact, it is homeomorphic to a Cantor set.
\end{itemize}
\end{lemma}

The proof is a folklore and it can be found in \cite{BFP13}.
Kakeya thought that unions of line segments, as in part (a), and the Cantor set, as in part (b), are the only two topological possibilities for the set of all subsums \eqref{eq:achievement} (a.k.a.\@ as the \emph{achievement set}) of a positive series $\sum_{n}x_n$.
Alek Weinstein and Boris Shapiro \cite{WS80} might have been the first to notice that there is a third possibility, which has been called the \emph{Cantorval}. Today it is known that the list of possibilities is complete \cite{GN88,NS00}.

Also observe that the condition \eqref{eq:Kak2}, satisfied for all large $n$, says nothing about the Lebesgue measure (or the Hausdorff dimension) of \eqref{eq:achievement}.
For instance, $x_n=1/3^n$ leads to the triadic Cantor set of measure $0$, while $x_n=1/(2^n-1)$ leads to a so-called \emph{fat Cantor set}, the one of positive measure.

We can easily prove a minor (but very useful) generalization of Kakeya's result.

\begin{lemma}\label{lm:new}
Let $X_1,X_2,X_3,\ldots$ be finite subsets of $[0,\infty)$, each with at least two elements, such that the series $\sum_n \max X_n$ converges. For every $n\in\N$ let $\Delta_n$ and $\delta_n$ respectively be the largest and the smallest length of intervals into which $X_n$ subdivides $[\min X_n, \max X_n]$. Also denote
\[ r_n := \sum_{k=n+1}^{\infty} \big( \max X_k - \min X_k \big). \]
\begin{itemize}
\item[(a)] If
\begin{equation}\label{eq:cond1}
r_n\geq \Delta_n
\end{equation}
holds for every sufficiently large $n\in\N$, then
\begin{equation}\label{eq:sumset}
\Big\{ \sum_{n=1}^{\infty} x_n \,:\, x_n\in X_n \text{ for every } n\in\N \Big\}
\end{equation}
is a finite union of nondegenerate bounded closed intervals. 
Moreover, if condition \eqref{eq:cond1} holds for every $n\in\N$, then the set \eqref{eq:sumset} is a single interval, namely 
\begin{equation}\label{eq:theinterval}
\Big[\sum_{n=1}^{\infty}\min X_n,\sum_{n=1}^{\infty}\max X_n\Big]. 
\end{equation}
\item[(b)] If
\begin{equation}\label{eq:cond2}
r_n < \delta_n
\end{equation}
holds for every sufficiently large $n\in\N$, then \eqref{eq:sumset} is a closed set with empty interior.
\end{itemize}
\end{lemma}

Note that Lemma \ref{lm:Kakeya} is recovered by taking $X_n=\{0,x_n\}$ for every $n\in\N$.
Variants of Lemma \ref{lm:new} have already been formulated and applied in connection with irrationality problems, cf.\@ \cite[Lemmae 14 \& 16]{KT25}, but here we tried to be even slightly more general. We also decided to include the detailed proof of Lemma \ref{lm:new}, for completeness. 

\begin{proof}
Denote
\[ a_n := \sum_{k=n+1}^{\infty} \min X_k, \quad b_n := \sum_{k=n+1}^{\infty} \max X_k \]
for each $n\in\N\cup\{0\}$.

(a) Let us first suppose that \eqref{eq:cond1} holds for every index $n\in\N$. Take some number $x$ from the interval \eqref{eq:theinterval}.
Let us inductively construct $x_n\in X_n$ for $n=1,2,3,\ldots$ such that for every $N\in\N\cup\{0\}$ one has
\begin{equation}\label{eq:proofnew}
x - \sum_{n=1}^{N} x_n \in [a_N,b_N].
\end{equation}
The claim is trivial for $N=0$, by the choice of $x$ and since the interval $[a_0,b_0]$ is precisely \eqref{eq:theinterval}.
Now take some $N\in\N\cup\{0\}$ and suppose that $x_n\in X_n$, $n=1,2,\ldots,N$, have already been chosen to satisfy \eqref{eq:proofnew}.
Condition \eqref{eq:cond1} can be written as $\Delta_{N+1} \leq b_{N+1} - a_{N+1}$ and it guarantees
\[ [a_N,b_N] \subseteq \bigcup_{y\in X_{N+1}} \big(y + [a_{N+1},b_{N+1}]\big), \]
so we can define $x_{N+1}$ to be the smallest element of $X_{N+1}$ such that 
\[ x - \sum_{n=1}^{N} x_n - x_{N+1} \in [a_{N+1},b_{N+1}]. \]
This is clearly just \eqref{eq:proofnew} again, with $N+1$ in place of $N$, so the construction step is complete.
From the convergence of $\sum_n \max X_n$ we know 
\[ \lim_{N\to\infty}a_N=0=\lim_{N\to\infty}b_N, \] 
so \eqref{eq:proofnew} and the sandwich theorem imply
\[ x = \sum_{n=1}^{\infty} x_n, \]
as desired.

Now assume that there only exists an index $m\in\N$ such that \eqref{eq:cond1} holds for every $n>m$. 
We can write \eqref{eq:sumset} as the sumset
\begin{align*}
& \Big\{ \sum_{n=1}^{m} x_n \,:\, x_n\in X_n \text{ for every } 1\leq n\leq m \Big\} \\
& + \Big\{ \sum_{n=m+1}^{\infty} x_n \,:\, x_n\in X_n \text{ for every } n>m \Big\}. 
\end{align*}
The first summand above is a finite set, while the second one equals $[a_m,b_m]$, from the previously established case applied to the sets $X_{m+1},X_{m+2},\ldots$.
Thus, \eqref{eq:sumset} is really a union of finitely many translates of the same nondegenerate line segment.

(b) Let us again begin by assuming that the relevant condition, which is now \eqref{eq:cond2}, holds for every $n\in\N$.
It is clear that the set \eqref{eq:sumset} is, for every $N\in\N$, contained in the union of finitely many intervals,
\[ I_N(x_1,\ldots,x_N) := \sum_{n=1}^{N} x_n + [a_N,b_N], \]
defined as the $N$-tuple $(x_1,\ldots,x_N)$ ranges over $X_1\times\cdots\times X_N$.
We claim that all these intervals are mutually disjoint. Namely, if we take two different $N$-tuples, $(x_1,\ldots,x_N)$ and $(x'_1,\ldots,x'_N)$, and find an index $1\leq l\leq N$ such that (without loss of generality)
\[ x_1=x'_1, \quad \ldots,\quad x_{l-1}=x'_{l-1},\quad x_l<x'_l, \]
then
\[ I_N(x_1,\ldots,x_N) \subseteq \sum_{n=1}^{l-1} x_n + x_l + [0,b_l], \quad I_N(x'_1,\ldots,x'_N) \subseteq \sum_{n=1}^{l-1} x_n + x'_l + [a_l,\infty). \]
These two are disjoint because
\[ x'_l - x_l \geq \delta_l > r_l = b_l - a_l, \]
i.e., $x_l + b_l < x'_l + a_l$.
Now we know that every interval contained in \eqref{eq:sumset} has length at most $b_N-a_N=r_N$. The claim about the interior follows since $N$ was arbitrary and $\lim_{N\to\infty}r_N=0$.
Also note that \eqref{eq:sumset} is, in fact, equal to
\begin{equation}\label{eq:alternativerep}
\bigcap_{N=1}^{\infty} \bigcup_{(x_1,\ldots,x_N)\in X_1\times\cdots\times X_N} I_N(x_1,\ldots,x_N) 
\end{equation} 
and consequently it is a closed (and even compact) set.

Finally, if there is an index $m\in\N$ such that \eqref{eq:cond2} holds for $n>m$, then
\[ \Big\{ \sum_{n=m+1}^{\infty} x_n \,:\, x_n\in X_n \text{ for every } n>m \Big\} \]
is a closed set with empty interior, by the previous case applied to $X_{m+1},X_{m+2},\ldots$.
Then \eqref{eq:sumset} is a finite union of closed sets with empty interiors, so the whole set itself still has empty interior by the Baire category theorem.
\end{proof}

We can describe part (a) of Lemma \ref{lm:new} informally as follows. In some sense the sets $X_1,X_2,X_3,\ldots$ give sufficiently many alternatives for the series terms $x_1,x_2,x_3,\ldots$ and the interior starts building up as soon as these sets are sufficiently dense. 

\begin{remark}\label{rem:stronger}
Despite being a somewhat basic result, it is possible to apply Lemma \ref{lm:new} in novel ways in concrete situations like ours. The result is useful and nontrivial already when we replace \eqref{eq:cond1} with a stronger condition:
\[ \max X_{n+1} - \min X_{n+1} \geq \Delta_n \]
for all indices $n\in\N$.
\end{remark}

\begin{remark}\label{rem:measure}
If \eqref{eq:cond2} is satisfied for every index $n\in\N$, then the Lebesgue measure of the set \eqref{eq:sumset} can be read off from the previous proof. Namely, after writing it as \eqref{eq:alternativerep}, we see that its measure is equal to
\[ \lim_{N\to\infty} \sum_{(x_1,\ldots,x_N)\in X_1\times\cdots\times X_N} |I_N(x_1,\ldots,x_N)|
= \lim_{N\to\infty} |X_1|\cdots |X_N| \,r_N. \]
\end{remark}


\section{Proof of Theorem \ref{thm:main}}

Fix arbitrary numbers $0<\theta<M$; it is sufficient to show that the set \eqref{eq:thewholeset} contains to whole line segment $[\theta,M]$.
Partition the set of positive integers as $\N=\cup_{j=1}^{\infty}S_j$, where
\[ S_j := 2^{j-1} (2\N-1) = \{ 2^{j-1} \cdot 1,\, 2^{j-1} \cdot 3,\, 2^{j-1} \cdot 5,\, \ldots \}. \]
The idea is to apply Lemma \ref{lm:new} to the sets
\begin{equation}\label{eq:defofX}
X_j := \Big\{ \sum_{n\in S_j} \frac{1}{\prod_{i=1}^{f(n)}(n+i)} \,:\, f\in\mathcal{F}_j \Big\}, \quad j=1,2,3,\ldots,
\end{equation}
where $\mathcal{F}_j$ is a carefully chosen finite collection of functions from $S_j$ to $\N$.
The construction of $\mathcal{F}_j$ will in turn use Lemma \ref{lm:Kakeya} for a series with terms indexed by $S_j$, so our argument can be thought of as a two-step approximation of a real number from $[\theta,M]$.

Let us first observe that for every $j\in\N$, $j\geq2$ the series obtained from 
\[ \sum_{n\in S_j} \frac{1}{\prod_{i=1}^{j}(n+i)} \]
by summing over $S_j$ in the ascending order, satisfies the condition \eqref{eq:Kak1} from Lemma \ref{lm:Kakeya} for all sufficiently large indices. Namely, for every $l\in\N$ we have
\[ \sum_{k=l+1}^{\infty} \frac{1}{\prod_{i=1}^{j}(2^{j-1}(2k-1)+i)} 
\geq \sum_{k=l+1}^{\infty} \frac{1}{(2^j k)^j} 
\geq \frac{1}{2^{j^2}} \int_{l+1}^{\infty} \frac{\textup{d}t}{t^j} 
= \Omega(l^{-j+1}), \]
which is greater than 
\[ \frac{1}{\prod_{i=1}^{j}(2^{j-1}(2l-1)+i)} = O(l^{-j}) \]
as soon as $l$ is sufficiently large (depending on $j$).
By Lemma \ref{lm:Kakeya} we know that the set of its subsums 
\[ \Big\{ \sum_{n\in T} \frac{1}{\prod_{i=1}^{j}(n+i)} \,:\, T\subseteq S_j \Big\} \]
is a finite union of nondegenerate line segments. One of these segments contains $0$ and denote its length by $\varepsilon'_j>0$.
Finally, take 
\[ \varepsilon_j := \min\Big\{ \varepsilon'_j,\frac{\theta}{2^j} \Big\} \] 
for each $j\geq 2$.

We begin by the construction of $\mathcal{F}_1$.
Take the smallest positive integer $m_1$ such that 
\[ \vartheta_1 := \frac{M}{m_1} < \frac{\varepsilon_2}{4} \]
and partition the segment $[\theta/2, M + \theta/2]$ equidistantly using $m_1+1$ points
\[ p_{1,k} := \frac{\theta}{2} + k \vartheta_1, \quad k=0,1,\ldots,m_1. \]
The series $\sum_{n\in S_1} 1/(n+1)$ is a variant of the harmonic series and it is well-known to have subsums
\[ \Big\{ \sum_{n\in T} \frac{1}{n+1} \,:\, T\subseteq S_1 \Big\} \]
equal to any prescribed positive number; see for instance \cite[\S 2]{Kov24}.
Thus, for every $0\leq k\leq m_1$ we can find a set $T'_{1,k}\subseteq S_1$ such that 
\[ \sum_{n\in T'_{1,k}} \frac{1}{n+1} = p_{1,k} \]
and then a finite set $T_{1,k}\subseteq T'_{1,k}$ such that
\[ \sum_{n\in T'_{1,k}\setminus T_{1,k}} \frac{1}{n+1} < \frac{\vartheta_1}{8}. \]
For each $0\leq k\leq m_1$ the function $f_{1,k}\colon S_1\to\N$ is defined to be
\begin{equation}\label{eq:growth1}
f_{1,k}(n) := 1 \quad \text{for every } n\in T_{1,k},
\end{equation}
while its values $f_{1,k}(n)$, $n\in S_1\setminus T_{1,k}$, can be quite arbitrary, and we only take care that
\begin{equation}\label{eq:growth2}
f_{1,k}(n) \geq n+1 \quad \text{for every } n\in S_1\setminus T_{1,k}
\end{equation}
and
\[ \sum_{n\in S_1\setminus T_{1,k}} \frac{1}{\prod_{i=1}^{f_{1,k}(n)}(n+i)} < \frac{\vartheta_1}{8}. \]
That way
\[ \Big| \sum_{n\in S_1} \frac{1}{\prod_{i=1}^{f_{1,k}(n)}(n+i)} - p_{1,k} \big| < \frac{\vartheta_1}{4} \quad \text{for } k=0,1,\ldots,m_1. \]
Setting
\[ \mathcal{F}_1 := \big\{f_{1,k} : k\in\{0,1,\ldots,m_1\} \big\} \]
and defining $X_1$ as in \eqref{eq:defofX}, we conclude (in the notation of Lemma \ref{lm:new}):
\begin{equation}\label{eq:ineqq1} 
0 < \min X_1 \leq \frac{\theta}{2} + \frac{\vartheta_1}{4}, \quad
\max X_1 \geq M + \frac{\theta}{2} - \frac{\vartheta_1}{4}, \quad
\Delta_1 < 2\vartheta_1. 
\end{equation}

Now we turn to the construction of $\mathcal{F}_j$ for any fixed positive integer $j\geq2$.
Take the smallest positive integer $m_j$ such that 
\[ \vartheta_j := \frac{\varepsilon_j}{m_j} < \frac{\varepsilon_{j+1}}{4} \]
and partition the segment $[0,\varepsilon_j]$ equidistantly using $m_j+1$ points
\[ p_{j,k} := k \vartheta_j, \quad k=0,1,\ldots,m_j. \]
We already know that for every $0\leq k\leq m_j$ there exists a set $T'_{j,k}\subseteq S_j$ such that 
\[ \sum_{n\in T'_{j,k}} \frac{1}{\prod_{i=1}^{j}(n+i)} = p_{j,k}. \]
Then we choose a finite set $T_{j,k}\subseteq T'_{j,k}$ such that
\[ \sum_{n\in T'_{j,k}\setminus T_{j,k}} \frac{1}{\prod_{i=1}^{j}(n+i)} < \frac{\vartheta_j}{8} \]
and define a function $f_{j,k}\colon S_j\to\N$ by
\begin{equation}\label{eq:growth3}
f_{j,k}(n) := j \quad \text{for every } n\in T_{j,k}, 
\end{equation}
only requiring
\begin{equation}\label{eq:growth4}
f_{j,k}(n) \geq n+j \quad \text{for every } n\in S_j\setminus T_{j,k} 
\end{equation}
and
\[ \sum_{n\in S_j\setminus T_{j,k}} \frac{1}{\prod_{i=1}^{f_{j,k}(n)}(n+i)} < \frac{\vartheta_j}{8}. \]
Consequently,
\[ \Big| \sum_{n\in S_j} \frac{1}{\prod_{i=1}^{f_{j,k}(n)}(n+i)} - p_{j,k} \big| < \frac{\vartheta_j}{4} \quad \text{for } k=0,1,\ldots,m_j. \]
After setting
\[ \mathcal{F}_j := \big\{f_{j,k} : k\in\{0,1,\ldots,m_j\} \big\} \]
and again defining $X_j$ as in \eqref{eq:defofX}, we can write:
\begin{equation}\label{eq:ineqq2} 
0 < \min X_j \leq \frac{\vartheta_j}{4}, \quad
\max X_j \geq \varepsilon_j - \frac{\vartheta_j}{4}, \quad
\Delta_j < 2\vartheta_j. 
\end{equation}

Note that the conditions of Lemma \ref{lm:new} are fulfilled, even in the stronger form from Remark \ref{rem:stronger}:
\[ \max X_{j+1} - \min X_{j+1} \stackrel{\eqref{eq:ineqq2}}{\geq} \varepsilon_{j+1} - \frac{\vartheta_{j+1}}{2} 
\geq \frac{\varepsilon_{j+1}}{2} > 2\vartheta_j \stackrel{\eqref{eq:ineqq1},\eqref{eq:ineqq2}}{>} \Delta_j \]
for every $j\in\N$. Also note that 
\[ \sum_{j=1}^{\infty} \min X_j \stackrel{\eqref{eq:ineqq1},\eqref{eq:ineqq2}}{\leq} \frac{\theta}{2} + \frac{1}{4} \sum_{j=1}^{\infty} \vartheta_j 
\leq \frac{\theta}{2} + \frac{1}{4} \sum_{j=1}^{\infty} \frac{\theta}{2^j} < \theta \]
and
\[ \sum_{j=1}^{\infty} \max X_j \stackrel{\eqref{eq:ineqq1},\eqref{eq:ineqq2}}{\geq} M + \frac{\theta}{2} + \sum_{j=2}^{\infty} \varepsilon_j  - \frac{1}{4} \sum_{j=1}^{\infty} \vartheta_j 
> M + \frac{\theta}{2} - \frac{1}{4} \sum_{j=1}^{\infty} \frac{\theta}{2^j} > M, \]
so the interval \eqref{eq:theinterval} certainly contains $[\theta,M]$.
It remains to apply Lemma \ref{lm:new} to the sets $X_1,X_2,X_3,\ldots$. For any given $x\in[\theta,M]$ it produces a sequence of functions $f_j\in\mathcal{F}_j$ such that
\[ x = \sum_{j=1}^{\infty} \sum_{n\in S_j} \frac{1}{\prod_{i=1}^{f_j(n)}(n+i)}. \]
We then have
\[ x = \sum_{n\in\N} \frac{1}{\prod_{i=1}^{f(n)}(n+i)}, \]
once we define $f\colon\N\to\N$ as $f(n) := f_j(n)$ for $j\in\N$ and $n\in S_j$.
Note that, for every $N\in\N$, the functions $f_1,\ldots,f_N$ attain finitely many values that are at most $N$, while the functions $f_{N+1},f_{N+2},\ldots$ only attain values greater than $N$; all thanks to \eqref{eq:growth1}, \eqref{eq:growth2}, \eqref{eq:growth3}, and \eqref{eq:growth4}. This clearly shows $\lim_{n\to\infty}f(n)=\infty$.


\section{Proof of Theorem \ref{thm:cantor}}

Since we want to prove that the set \eqref{eq:increasingset} has zero measure, and not just empty interior, we do not apply Lemma \ref{lm:new} directly, but rather proceed similarly as in its proof (cf.\@ Remark \ref{rem:measure}).

For any fixed $M\in\N$ we only consider the intersection of \eqref{eq:increasingset} with $[0,M]$ and show that this intersection has measure $0$. Take $K\in\N$ so that
\begin{equation}\label{eq:choiceK}
\sum_{n=1}^{K} \frac{1}{n+1} > M,
\end{equation}
which exists due to divergence of the harmonic series.
Fix some positive integer $N$ such that $N\geq K$.

For every $m\in\N$ and every tuple $(t_1,\ldots,t_{m-1})$ from
\[ \mathcal{T}_{N,m} := \{ (t_1,\ldots,t_{m-1})\in\N^{m-1} \,:\, t_1\leq t_2\leq \cdots \leq t_{m-1}\leq N \} \]
we define the interval
\[ I_{N,m}(t_1,\ldots,t_{m-1}) := \sum_{n=1}^{m-1} \frac{1}{\prod_{i=1}^{t_i}(n+i)} + \Big[ 0, \sum_{n=m}^{\infty} \frac{1}{\prod_{i=1}^{N+1}(n+i)} \Big]. \]
We also define
\[ I_{N,1} := \Big[ 0, \sum_{n=1}^{\infty} \frac{1}{\prod_{i=1}^{N+1}(n+i)} \Big] \]
and later interpret $I_{N,1}()$, with an empty list of parameters $t_1,t_2,\ldots$, precisely as $I_{N,1}$.
To every increasing function $f\colon\N\to\N$ satisfying $\lim_{n\to\infty}f(n)=\infty$ we can associate unique $m\in\N$ and $(t_1,\ldots,t_{m-1})\in\mathcal{T}_{N,m}$ such that
\[ \underbrace{f(1)}_{=t_1}\leq \underbrace{f(2)}_{=t_2}\leq \cdots \leq \underbrace{f(m-1)}_{=t_{m-1}}\leq N < f(m) \leq f(m+1)\leq \cdots, \]
so \eqref{eq:increasingset} intersected with $[0,M]$ is contained in
\begin{equation}\label{eq:iunion}
\bigcup_{m\in\N} \bigcup_{(t_1,\ldots,t_{m-1})\in\mathcal{T}_{N,m}} I_{N,m}(t_1,\ldots,t_{m-1}) \cap [0,M]. 
\end{equation} 
Our goal is to show that the measures of \eqref{eq:iunion} tend to $0$ as $N\to\infty$.

The length of the interval $I_{N,m}(t_1,\ldots,t_{m-1})$ (including $I_{N,1}$) only depends on $N$ and $m$ and it equals
\[ \sum_{n=m}^{\infty} \frac{1}{\prod_{i=1}^{N+1}(n+i)}
= \frac{1}{N} \sum_{n=m}^{\infty} \Big( \frac{1}{\prod_{i=1}^{N}(n+i)} - \frac{1}{\prod_{i=1}^{N}(n+1+i)} \Big)
= \frac{1}{N} \frac{1}{\prod_{i=1}^{N}(m+i)}. \]
Cardinality of the collection $\mathcal{T}_{N,m}$ is the number of $(m-1)$-combinations with repetitions from a set of size $N$, i.e., 
\[ |\mathcal{T}_{N,m}| = \binom{N+(m-1)-1}{m-1} = \binom{m+N-2}{N-1}. \]
Consequently, for every $m\in\N$ such that $1\leq m\leq N$ we get the bound
\begin{align}
& \sum_{m=1}^{N} \sum_{(t_1,\ldots,t_{m-1})\in\mathcal{T}_{N,m}} |I_{N,m}(t_1,\ldots,t_{m-1})| \nonumber \\
& = \sum_{m=1}^{N} \frac{1}{(m+N)N!}\frac{(m+N-2)(m+N-3)\cdots m}{(m+N-1)(m+N-2)\cdots (m+1)} < \frac{1}{N!}. \label{eq:ttt1}
\end{align}

For those $m\in\N$ that are greater than $N$ we need to be more economical with estimation, so observe that the tuples $(t_1,\ldots,t_{m-1})$ from $\mathcal{T}_{N,m}$ that begin as
\[ t_1 = t_2 = \cdots = t_K = 1 \]
push the interval $I_{N,m}(t_1,\ldots,t_{m-1})$ right of the number $M$, due to $m-1\geq N\geq K$ and the choice \eqref{eq:choiceK} of the parameter $K$.
Thus, for a fixed $m>N$, at most
\[ \binom{m+N-2}{N-1} - \binom{m-K+N-2}{N-1} \]
intervals in \eqref{eq:iunion} have nonempty intersections with $[0,M]$.
The last difference of binomials is easily seen to be at most
\[ \frac{(N+K-2)(m+N-2)^{N-2}}{(N-2)!}, \]
so 
\begin{align}
& \sum_{m=N+1}^{\infty} \sum_{(t_1,\ldots,t_{m-1})\in\mathcal{T}_{N,m}} |I_{N,m}(t_1,\ldots,t_{m-1}) \cap [0,M]| \nonumber \\
& \leq \sum_{m=N+1}^{\infty} \frac{2N(2m)^{N-2}}{(N-2)!} \frac{1}{Nm^N} 
= \frac{2^{N-1}}{(N-2)!} \sum_{m=N+1}^{\infty}\frac{1}{m^2} < \frac{2^N}{(N-2)!}. \label{eq:ttt2}
\end{align}

Finally, we conclude from \eqref{eq:ttt1} and \eqref{eq:ttt2} that the measure of \eqref{eq:iunion} is at most
\[ \frac{1}{N!} + \frac{2^N}{(N-2)!}, \]
which clearly converges to $0$ as $N\to\infty$.
Our set of interest \eqref{eq:increasingset} is contained in the union over $M\geq1$ of the intersections over $N\geq K$ of the sets \eqref{eq:iunion}. We conclude that \eqref{eq:increasingset} is a measurable set of Lebesgue measure zero.


\section*{Acknowledgments}
This work was supported by the Croatian Science Foundation under the project number HRZZ-IP-2022-10-5116 (FANAP).
The authors are grateful to Thomas Bloom for founding and managing the website \cite{EP}.


\bibliography{RationalErdosGraham}{}
\bibliographystyle{plainurl}

\end{document}